\documentclass[12pt]{article}
\usepackage{amsmath}
\usepackage{amscd}
\usepackage{amsthm}
\usepackage{amssymb}
\usepackage[latin1]{inputenc}
\usepackage[all]{xy}
\def\email#1{{\tt #1}}
\def\SetSize{\fontsize{12}{14.4}\selectfont}

\newtheoremstyle{theorem}
  {}
  {}
  {\itshape}
  {}
  {}
  {.}
  {.5em}
  {}

\newtheoremstyle{definition}
  {}
  {}
  {}
  {}
  {\MakeUppercase}
  {.}
  {.5em}
  {}

\theoremstyle{theorem}
\newtheorem{theorem}{THEOREM}[section]

\newtheorem{corollary}[theorem]{COROLLARY}
\newtheorem{proposition}[theorem]{PROPOSITION}

\theoremstyle{definition}
\newtheorem{remark}[theorem]{REMARK}

\newtheorem{exercise*}[enumiv]{EXERCISE}

\let\frak=\mathfrak

\let\phi=\varphi

\def\Cok{\operatorname{Cok}}

\def\grade{\operatorname{grade}}

\def\Ext{\operatorname{Ext}}
\def\rank{\operatorname{rank}}

\def\Ker{\operatorname{Ker}}
\def\Coker{\operatorname{Coker}}

\def\Im{\operatorname{Im}}

\let\oldbigwedge\bigwedge
\def\BIGwedge{{\textstyle\oldbigwedge}}
\def\medwedge{{\scriptstyle\oldbigwedge}}
\def\bigwedge{\mathchoice{\BIGwedge}{\BIGwedge}{\medwedge}{}}

\let\iso=\cong

\let\epsilon=\varepsilon
\let\tilde=\widetilde

\begin{document}
\thispagestyle{empty} \vspace*{1.5in}
{\fontsize{14}{16.8}\selectfont \noindent \uppercase{Length
Formulas for the homology of generalized Koszul complexes }\par}

\SetSize \vspace{2\baselineskip}

\noindent\uppercase{Bogdan Ichim}, Universit\"at Osnabr\"uck, FB
Mathematik/Informatik, 49069 Osna\-br\"uck, Germany,
\email{bogdan.ichim@mathematik.uni-osnabrueck.de}\\
Institute of Mathematics, C.P. 1-764, 70700 Bucharest, Romania,\\
\email{bogdan.ichim@imar.ro}
\\[1\baselineskip]
\noindent\uppercase{Udo Vetter}, Universit\"at Oldenburg, Institut
für Mathematik, 26111 Oldenburg, Germany,
\email{vetter@\allowbreak mathematik.uni-oldenburg.de}

\vspace{4\baselineskip plus 1 \baselineskip minus 1\baselineskip}

\begin{abstract}
Let $M$ be a finite module over a noetherian ring $R$ with a free
resolution of length $1$.
 We consider the generalized Koszul
complexes $\mathcal{C}_{\bar\lambda}(t)$ associated with a map
$\bar\lambda:M\to\mathcal{H}$ into a finite free $R$-module
$\mathcal{H}$ (see [IV], section 3), and investigate the homology
of $\mathcal{C}_{\bar\lambda}(t)$ in the special setup when
$\grade I_M=\rank M=\dim R$. ($I_M$ is the first non-vanishing
Fitting ideal of $M$.) In this case the (interesting) homology of
$\mathcal{C}_{\bar\lambda}(t)$ has finite length, and we deduce
some length formulas. As an application we give a short algebraic
proof of an old theorem due to Greuel (see [G], Proposition 2.5).
We refer to [HM] where one can find another proof by similar
methods.

\end{abstract}

\section*{Introduction}

Let $R$ be a noetherian ring and $M$ an $R$-module with a
presentation
$$
\CD 0@>>>\mathcal{F}@>\chi>>\mathcal{G}@>>>M@>>>0
\endCD
$$
where $\mathcal{F}$, $\mathcal{G}$ are free modules of finite
rank. We consider an $R$-homomorphism
$\bar\lambda:M\to\mathcal{H}$ into a finite free $R$-module
$\mathcal{H}$. With $\bar\lambda$ we associate the generalized
Koszul complexes $\mathcal{C}_{\bar\lambda}(t)$. The main object
of our interest is the homology of these complexes in the special
setup when $\grade I_\chi=\rank M=\dim R$. Here $I_\chi$ denotes
the ideal generated by the maximal minors of a matrix representing
$\chi$.

In section 1 we give a short survey of generalized Koszul
complexes and Koszul bicomplexes. Section 2 contains some results
concerning the grade sensitivity of these and some related
complexes.

The main part of the paper is concentrated in section 3. In case
$\grade I_\chi$ has the greatest possible value $\rank M+1$,
Theorem 6.10 in [IV] provides a satisfactory description of the
homology of $\mathcal C_{\bar{\lambda}}(t)$ depending on the grade
$h$ of $I_\lambda$, $\lambda:\mathcal G\to \mathcal H$ being the
corresponding lifted map. In a sense the situation described above
is another extremal case. Here the homology modules $\tilde H^i$
of $\mathcal C_{\bar{\lambda}}(t)$ have finite length for
$i\le\min(h-1,2\rank M)$, and they are connected by some length
formulas which mainly involve the symmetric powers of the cokernel
of $\chi^*$ (Theorem \ref{SubMaximalCaseHomology} and Corollary
\ref{length formulas}).

In section 4 we treat the very special situation of a
quasi-homogeneous isolated complete intersection singularity. We
can extend the length formulas of section 3, in order to give a
purely algebraic proof of an old theorem due to Greuel (see [G],
2.5). The formulas and -- at least implicitly -- an algebraic
proof of Greuel's Theorem are also contained in [HM].

We use [E] for a general reference in commutative algebra. As
there we denote by $\bigwedge M$ ($\bigwedge^p M$) the exterior
power algebra ($p$th exterior power) of an $R$-module $M$, by
$S(M)$ ($S_p(M)$) the symmetric power algebra ($p$th symmetric
power) of $M$, and by $D(M)$ ($D_p(M)$) the divided power algebra
($p$th divided power) of $M$. The reader should always have in
mind that for a finite free $R$-module $F$ there are canonical
isomorphisms $D(F^*)\iso S(F)^*$ and $D(F)\iso S(F^*)^*$ which we
shall use implicitly several times in the following. Here
$S(F)^*=\oplus S_p(F)^*$ is the so called graded dual of $S(F)$.

\section{Koszul Complexes and Koszul Bicomplexes}\label{Koszul
Complexes}

Generalizations of the classical Koszul complex have been known
since a long time, e.\ g. the complexes due to Eagon and Northcott
or those introduced by Buchsbaum and Rim. We refer to [E], Chapter
A2.6, for a comprehensive treatment. Some related material may
also be found in [BV1] and [BV2].

So the following constructions, generalizing the classical Koszul
complex and its dual, are not really new. They should be viewed as
appropriate components of the Koszul bicomplex, we define at the
end of this section and which is the main tool of our
investigations.

Let $R$ be a commutative ring and $\psi: G\to F $ a homomorphism
of an $R$-module $G$ into a finite free $R$-module $F$. Let
$f_1,\ldots,f_m$ be a basis of $F$ and $f_1^*,\ldots,f_m^*$ the
dual basis of $F^*$. Then we may consider $\psi$ as an element of
$\bigwedge G^*\otimes S(F)$ with the presentation $
\psi=\sum_{j=1}^m\psi^*(f_j^*)\otimes f_j. $ We set
\begin{equation*}
\partial_\psi(y_1\wedge\ldots\wedge y_n\otimes z) =
\sum_j\big(y_1\wedge\ldots\wedge
y_n\leftharpoonup\psi^*(f_j^*)\big)\otimes z\cdot
f_j.\end{equation*} for all $y_i\in G$ and $z\in S(F)^*$. The
right multiplication $\leftharpoonup$ of $\bigwedge G^*$ on
$\bigwedge G$ is given by
\begin{equation*}
y_1\wedge\ldots\wedge y_n\leftharpoonup y^*_1\wedge\ldots\wedge
y^*_p=\sum_{\sigma}\epsilon(\sigma) \det_{1\le i,j\le p}
(y^*_j(y_{\sigma(i)}))y_{\sigma(p+1)}\wedge\ldots\wedge
y_{\sigma(n)}
\end{equation*}
for $y_1,\ldots ,y_n\in G$ and $y^*_1,\ldots ,y^*_p\in G^*$, where
$\sigma$ runs through the set $\mathfrak{S}_{n,p}$ of permutations
of $n$ elements which are increasing on the intervals $[1,p]$ and
$[p+1,n]$. The {\it generalized Koszul complexes} we have in mind,
are the complexes
\begin{align*}
\mathcal C_\psi(t):\quad \cdots\to\bigwedge^{t+m+p} G \otimes
S_{p}(F)^*\stackrel{\partial_\psi}{\to}\cdots
\stackrel{\partial_\psi}{\to}\bigwedge^{t+m} G & \otimes
S_{0}(F)^* \stackrel{\nu_{\psi}}{\to}\bigwedge^{t} G \otimes
S_{0}(F)\stackrel{\partial_\psi}{\to}\\ &\cdots
\stackrel{\partial_\psi}{\to} \bigwedge^{0} G \otimes S_{t}(F) \to
0;
\end{align*}
$\nu_{\psi}$ is the right multiplication by $\psi^*(f^*_1)\wedge
\ldots \wedge \psi^*(f^*_m)\in \bigwedge G^*$.
\smallskip

Similarly we can associate complexes with a map $\phi:H\to G$ from
a finite free $R$-module $H$ into $G$, generalizing the dual
version of the classical Koszul complex. Let $h_1,\ldots,h_l$ be a
basis of $H$ and $h_1^*,\ldots,h_l^*$ the dual basis of $H^*$.
Then $\phi$, as an element of $S(H^*)\otimes\bigwedge G$, has the
presentation $ \phi=\sum_{j=1}^l h_j^*\otimes \phi(h_j).$ Define
$d_\phi$ to be the left multiplication by $\phi$ on $D(H)\otimes
\bigwedge G$, i.e.
$$d_\phi(x_1^{(k_1)}\ldots x_p^{(k_p)}\otimes
y)=\sum_j x_1^{(k_1)}\ldots x_j^{(k_j-1)}\ldots x_p^{(k_p)}\otimes
\phi(x_j)\wedge y$$ for $x_i\in H$ and $y\in\bigwedge G$, and on
$S(H^*)\otimes\bigwedge G $ (in an obvious way). We obtain the
family of complexes
\begin{align*}
\mathcal D_\phi(t):\quad 0\to D_{t}(H)\otimes \bigwedge^{0}
G\stackrel{d_\phi}{\to}\cdots \stackrel{d_\phi}{\to}
D_{0}(H)\otimes \bigwedge^{t} G & \stackrel{\nu^{\phi}}{\to}
S_{0}(H^*)\otimes \bigwedge^{t+l} G\stackrel{d_\phi}{\to}\cdots \\
& \stackrel{d_\phi}{\to}S_{p}(H^*)\otimes \bigwedge^{t+l+p} G
\to\cdots
\end{align*}
where $\nu^{\phi}$ is the left multiplication by $\phi(h_1)\wedge
\ldots \wedge \phi(h_m)\in \bigwedge G$.
\smallskip

Now let $\psi\circ\phi=0$. Then there is a simple associativity
formula (see Proposition 2.1 in [IV]) concerning the right and
left multiplications from above, which allows us to assemble the
complexes $\mathcal C_\psi(t)$ and $\mathcal D_\phi(t)$ to the
Koszul bicomplexes $\mathcal{K}_{.,.}(t)$ \scriptsize
$$
\CD
\vdots &&\vdots &&\vdots &&\vdots\\
@VVV @VVV @VVV @VVV\\
\cdots H\otimes\bigwedge^{t+m}G\otimes
F^*@>>>\bigwedge^{t+m+1}G\otimes F^*@>\pm\nu^{\phi}>>
\bigwedge^{t+l+m+1}G\otimes F^*@>>>H^*\otimes\bigwedge^{t+l+m+2}G\otimes F^*\cdots\\
@VVV @V\partial_{\psi}VV @V\partial_{\psi}VV @VVV\\
\cdots
H\otimes\bigwedge^{t+m-1}G@>d_{\phi}>>\bigwedge^{t+m}G@>\pm\nu^{\phi}>>
\bigwedge^{t+l+m}G@>d_{\phi}>>H^*\otimes\bigwedge^{t+l+m+1}G\cdots\\
@V\pm\nu_\psi VV @V\pm\nu_\psi VV @V\pm\nu_\psi VV @V\pm\nu_\psi VV\\
\cdots H\otimes\bigwedge^{t-1}G @>d_{\phi}>>\bigwedge^{t}G
@>\pm\nu^{\phi}>>
\bigwedge^{t+l}G@>d_{\phi}>>H^*\otimes\bigwedge^{t+l+1}G\cdots\\
@VVV @V\partial_{\psi}VV @V\partial_{\psi}VV @VVV\\
\cdots H\otimes\bigwedge^{t-2}G\otimes
F@>>>\bigwedge^{t-1}G\otimes F@>\pm\nu^{\phi}>>
\bigwedge^{t+l-1}G\otimes F@>>>H^*\otimes\bigwedge^{t+l}G\otimes F\cdots\\
@VVV @VVV @VVV @VVV\\
\vdots &&\vdots &&\vdots &&\vdots\\
\endCD
$$
\normalsize The rows in the upper half arise from $\mathcal
D_\phi(t+m+j)$ tensored with $S_j(F)^*$, $j=0,1,\ldots$, while the
rows below are built from $\mathcal D_\phi(t-j)$ tensored with
$S_j(F)$, $j=0,1,\ldots$; we abbreviate $d_\phi\otimes 1_{S(F^*)}$
and $d_\phi\otimes 1_{S(F)}$ to $d_\phi$, and correspondingly
$\nu^\phi\otimes 1_{S(F^*)}$ and $\nu^\phi\otimes 1_{S(F)}$ to
$\nu^\phi$. The columns are obtained analogously: in western
direction we have to tensorize $D_i(H)$ with $\mathcal
C_\psi(t-i)$, $i=0,1,\ldots$, while going east we must tensorize
$S_i(H^*)$ with $\mathcal C_\psi(t+l+i)$, $i=0,1,\ldots$; as
before we shorten the complex maps to $\partial_\psi$ and
$\nu_\psi$. The signs of $\nu^\phi$ and $\nu_\psi$ are determined
by the associativity formula.
\smallskip

A detailed and more general treatment of the generalized Koszul
complexes from above and of the Koszul bicomplexes just defined
may be found in [I].

\section{Grade Sensitivity}

The following contains some preparing material for the main
results in the next section. In a sense it is a matter of
generalizing part of the considerations of section 5 in [IV].

Our general assumption throughout the rest of this section will be
that $R$ is {\it noetherian}, that {\it $H$, $G$ and $F$ are free
$R$-modules of finite ranks $l$, $n$ and $m$, and that
$$
\CD H@>\phi>>G@>\psi>>F
\endCD
$$
is a complex.}  Although much of what we will do, holds formally
for any $l$, $n$ and $m$, the applications will refer to the case
in which $n\ge m$ and $n\ge l$. So {\it we suppose that $r=n-m\ge
0$, $s=n-l\ge 0$.} By $I_\psi$ ($I_\phi$) we denote the ideals in
$R$ generated by the maximal minors of a matrix representing
$\psi$ ($\phi$). We set $g=\grade I_\psi$ and $h=\grade I_\phi$.

Since $G$ is finitely generated, the generalized Koszul complexes
$\mathcal C_{\psi}(t)$ and $\mathcal D_{\phi}(t)$  have only a
finite number of non-vanishing components. To identify the
homology, we fix their graduations as follows: position 0 is held
by the leftmost non-zero module. As $\mathcal C_{\psi}(t)$ and
$\mathcal D_{\phi}(t)$ are isomorphic to well-known
generalizations of the classical Koszul complex (see [IV], section
3), their homology behaves grade sensitively in the sense of the
following theorem (see Theorem 5.2 in [IV]).
\begin{theorem}\label{E-N HOMOLOGY}Set
$C=\Coker \psi$ and $D=\Coker \phi^*$. Furthermore let
$S_0(D)=R/I_\phi$, $S_{-1}(D)=\bigwedge^{s+1}\Coker\phi$,
$S_0(C)=R/I_\psi$ and $S_{-1}(C)=\bigwedge^{r+1}\Coker\psi^*$.
Then the following hold.
\begin{enumerate}
\item[\rm{(a)}]$H^i(\mathcal D_{\phi}(t))=0$ for $i<h$. Moreover,
if $t\le s+1$ and $\grade I_k(\phi)\ge n-k+1$ for all $k$ with
$l\ge k\ge 1$, then $\mathcal D_{\phi}(t)$ is a free resolution of
$S_{s-t}(D)$. (If $-1\le t\le s+1$, then it suffices to require
that $\grade I_\phi\ge s+1$.) \item[\rm{(b)}]$H^i(\mathcal
C_{\psi}(t))=0$ for $i<g$. Moreover, if $t\ge -1$ and $\grade
I_k(\psi)\ge n-k+1$ for all $k$ with $m\ge k\ge 1$, then $\mathcal
C_{\psi}(t)$ is a free resolution of $S_{t}(C)$. (If $-1\le t\le
r+1$, then it suffices to require that $\grade I_\psi\ge r+1$.)
\end{enumerate}
Finally, if $I_{\phi}=R$ ($I_{\psi}=R$), then all sequences
$\mathcal D_{\phi}(t)$ ($\mathcal C_{\psi}(t)$) are split exact.
\end{theorem}

By $\mathcal{C}_{.,.}(t)$ we shall denote the bicomplex which is
the lower part of the Koszul bicomplex $\mathcal{K}_{.,.}(t)$ (the
rows below the second row in the last diagram of the previous
section). In other words,
$$C^{0,0}=\mathcal C_{0,0}(t)\ =\ \begin{cases}
D_t(H)\otimes\bigwedge^{0}G\otimes S_0(F)\quad &\text{if}\quad
0\le t,\\
S_0(H^*)\otimes\bigwedge^{t+l}G\otimes S_0(F)\quad
&\text{if}\quad -l\le t< 0,\\
S_{-t-l}(H^*)\otimes\bigwedge^{0}G\otimes S_0(F)\quad
&\text{if}\quad t< -l.\end{cases}
$$
The row homology at $C^{p,q}=\mathcal C_{p,q}(t)$ is denoted by
$H_\phi^{p,q}$, the column homology by $H_\psi^{p,q}$. Thus
$H_\phi^{p,0}$ is the $p$-th homology module of $\mathcal
D_{\phi}(t)$.

Set $N^p=\Ker {(C^{p,0}\overset {\partial_\psi}\to C^{p,1})}$. The
canonical injections $N^p\to C^{p,0}$ yield a complex morphism
$$ \CD 0@>>> N^0@>>>N^1 @>>>\cdots && N^p @>\bar{d_\phi}>>
N^{p+1}&&\cdots\\ &&
\parallel && @VVV && @VVV @VVV \\ 0@>>> C^{0,0}@>>>
C^{1,0}@>>>\cdots &&C^{p,0}@>d_\phi >> C^{p+1,0}&&\cdots
\endCD
$$ where the maps $\bar {d_\phi}$ are induced by $d_\phi$. The homology of
the first row $\mathcal N(t)$ at $N^p$ is denoted by $\bar{H}^p$.
We shall now investigate this homology under the assumption $r\le
g$.

For this purpose we extend $\mathcal{C}_{.,.}(t)$ by $\mathcal
N(t)$ to the complex $\tilde{\mathcal{C}}_{.,.}(t)$. So
$C^{p,-1}=\tilde{\mathcal C}_{p,-1}(t)=N^p$. We record some facts
about the homology of $\tilde{\mathcal{C}}_{.,.}(t)$. To avoid new
symbols, the column homology at $C^{p,q}$ is again denoted by
$H_\psi^{p,q}$ (actually it differs from that of
$\mathcal{C}_{.,.}(t)$ only at $C^{p,0}$). By construction
$H_\psi^{p,q}=0$ for $q=-1,0$. Furthermore we draw from Theorem
\ref{E-N HOMOLOGY},(a) that
\begin{equation}
H_\phi^{p,q}=0\quad\text{ for $p<h$ and $q\ne-1$}.
\end{equation}
Let $r\le g$. Then in case $0<p\le t$, we get from Theorem
\ref{E-N HOMOLOGY},(b) that
\begin{equation}
H_\psi^{p,q}=\begin{cases} \phantom{D_{t-p}(H)}0&\quad \text{for}\
q<\min(p-1,r),\\ D_{t-p}(H)\otimes S_p(C)&\quad \text{for}\
q=p,\end{cases} \end{equation} and if $0\le t<p$, we obtain
\begin{equation}
H_\psi^{p,q}=\begin{cases}
\phantom{S_{p-t-1}(H^*)}0&\quad\text{for}\
q<\min(p+l-2,r),\\S_{p-t-1}(H^*)\otimes S_{p+l-1}(C)&\quad
\text{for}\ q=p+l-1,\end{cases}
\end{equation}
where $C=\Coker\psi$ as in Theorem 2.1.

For $q\ge -1$ we consider the $q$th row $\tilde C_{.,q}$ of
$\tilde C_{.,.}$ and its image complex $\partial_\psi (\tilde
C_{.,q})$ in $\tilde C_{.,q+1}$. We set $E^{p,q}=H^p(\partial_\psi
(\tilde C_{.,q-1}))$ for $q\ge 0$. There are exact sequences
\begin{equation}
H_\phi^{i-(j+1),j} \to E^{i-(j+1),j+1}\to E^{i-j,j} \to
H_\phi^{i-j,j}
\end{equation}
if $H_\psi^{i-(j+1),j}=H_\psi^{i-j,j}=0$, and because of (2) and
(3) this holds if
\begin{align*}
0<i-(j+1)\le t\ \text{and}\
&j<\min\big(i-(j+1)-1,r\big)\ \text{or}\\
0\le t<i-(j+1)\ \text{and}\ &j<\min\big(i-(j+1)+l-2,r\big).
\end{align*}

\begin{theorem}\label{FUNDAMENTAL}Let $t\ge 0$ be an integer.
Assume that $ 1\le r\le g.$ Then, with the notation introduced
above, $\bar{H}^i= 0$ for $i=0,\ldots,\min(2,h-1)$. Set $C=\Coker
\psi$.
\begin{enumerate}
\item[\rm{(a)}] For $i$ odd, $3\le i < \min(h-1,2r,2t+2)$, one has
a natural exact sequence $$ 0\to \bar{H}^i\to
D_{t-\frac{i-1}2}(H)\otimes S_{\frac{i-1}2} (C) \to
H_\psi^{\frac{i+1}2,\frac{i-1}2} \to \bar{H}^{i+1} \to 0. $$
\item[\rm{(b)}] Suppose that $l>1$.
\begin{enumerate}
\item[\rm{(i)}] If $3\le 2t+1<h$, then $\bar{H}^{2t+1}=D_{0}(H)
\otimes S_{t} (C)$. \item[\rm{(ii)}] $\bar{H}^{i}=0$ for $2t+2\le
i< \min(h,2t+l+1)$. \item[\rm{(iii)}] If $2t+l+1<h$, then
$\bar{H}^{2t+l+1}=H_\psi^{t+1,t+l-1}$.
\end{enumerate}
\item[\rm{(c)}] For $i-l$ even, $2t+l+2\le i < \min(h-1,2r-l+2)$,
one has a natural exact sequence
$$
 0\to \bar{H}^i\to S_{\frac{i-l}2-t-1}(H^*)\otimes
S_{\frac{i+l}2-1} (C) \to H_\psi^{\frac{i-l}2+1,\frac{i+l}2-1} \to
\bar{H}^{i+1} \to 0.
$$
\end{enumerate}
\end{theorem}
\begin{proof} For the proof of the first statement we refer to [IV], Theorem 5.3.
Let $i$ be an odd integer, $3\le i< \min(h-1,2r,2t+2)$. Using (4),
we obtain the ``southwest'' isomorphisms
$$ \bar H^i=E^{i,0}\iso
E^{i-1,1}\iso\dots\iso E^{\frac {i+1}2,\frac
{i-1}2}\quad\text{and}\quad \bar H^{i+1}=E^{i+1,0}\iso
E^{i,1}\iso\dots\iso E^{\frac {i+3}2,\frac {i-1}2}.
$$
We abbreviate
$\partial_\psi^{p,q}=(C^{p,q}\overset{\partial_\psi}\to
C^{p,q+1}$) and $d_\phi^{p,q}=(C^{p,q}\overset{d_\phi}\to
C^{p+1,q})$. The diagram

$$ \CD && 0 && 0 && 0 && 0 \\ && @VVV @VVV @VVV @VVV \\
@>>>\Im\partial_\psi^{\frac{i-1}2,\frac{i-3}2}
@>>>\Im\partial_\psi^{\frac{i+1}2,\frac{i-3}2}
@>>>\Im\partial_\psi^{\frac{i+3}2,\frac{i-3}2}
@>>>\Im\partial_\psi^{\frac{i+5}2,\frac{i-3}2}\\ && @VVV @VVV @VVV @VVV \\
@>>>\Im\partial_\psi^{\frac{i-1}2,\frac{i-3}2}
@>>>\Ker\partial_\psi^{\frac{i+1}2,\frac{i-1}2}
@>>>\Ker\partial_\psi^{\frac{i+3}2,\frac{i-1}2}
@>>>\Ker\partial_\psi^{\frac{i+5}2,\frac{i-1}2}\\ && @VVV @VVV
@VVV @VVV
\\ && 0@>>> H_\psi^{\frac{i+1}2,\frac{i-1}2}@>>> 0@>>> 0\\ && &&  @VVV &&
&& &&\\ && &&  0 && && &&
\endCD
$$ is induced by $\mathcal{C}_{.,.}(t)$ and has exact columns. Its row
homology at $\Im\partial_\psi^{\frac{i+1}2,\frac{i-3}2}$ is
$E^{\frac{i+1}2,\frac{i-1}2}=\bar{H}^i$, and at
$\Ker\partial_\psi^{\frac{i+1}2,\frac{i-1}2}$ it coincides with
$D_{t-\frac{i-1}2}(H)\otimes S_{\frac{i-1}2}(C)$ as the following
diagram with exact columns and exact middle row shows (we write
$D_{t-\frac{i-1}2}$ for $D_{t-\frac{i-1}2}(H)$ and $S_{\frac
{i-1}2}$ for $S_{\frac {i-1}2}(C)$): \small
$$ \CD && 0 && 0 && 0 && 0
\\ && @VVV @VVV @VVV@VVV \\ &&
\Im\partial_\psi^{\frac{i-1}2,\frac{i-3}2}@>>>
\Ker\partial_\psi^{\frac{i+1}2,\frac{i-1}2}@>>>
\Ker\partial_\psi^{\frac{i+3}2,\frac{i-1}2}@>>>
\Ker\partial_\psi^{\frac{i+5}2,\frac{i-1}2}\\ && @VVV @VVV @VVV @VVV \\
0@>>> C^{\frac{i-1}2,\frac{i-1}2}@>>>
C^{\frac{i+1}2,\frac{i-1}2}@>>> C^{\frac{i+3}2,\frac{i-1}2}@>>>
C^{\frac{i+5}2,\frac{i-1}2}\\ &&@VVV @VVV @VVV @VVV \\ 0
@>>>D_{t-\frac{i-1}2}\otimes S_{\frac {i-1}2}@>>> \Im
\partial_\psi^{\frac{i+1}2,\frac{i-1}2}@>d_\phi>>
\Im\partial_\psi^{\frac{i+3}2,\frac{i-1}2}@>>>
\Im\partial_\psi^{\frac{i+5}2,\frac{i-1}2} \\ && @VVV @VVV @VVV @VVV \\
&&0 && 0 && 0 && 0
\endCD
$$
\normalsize

 In this diagram the row homology at
$\Im\partial_\psi^{\frac{i+1}2,\frac{i-1}2}$ vanishes since
$d_\phi^{\frac{i+1}2,\frac{i+1}2}$ is injective. So in the
preceding diagram the row homology at
$\Ker\partial_\psi^{\frac{i+3}2,\frac{i-1}2}$ is zero. Altogether
we obtain an exact sequence
\begin{equation*}\tag{$*$} 0\to \bar{H}^i\to
D_{t-\frac{i-1}2}(H)\otimes S_{\frac {i-1}2}(C)\to
H_\psi^{\frac{i+1}2,\frac{i-1}2}\to E^{\frac{i+3}2,\frac{i-1}2}\to
0.
\end{equation*}
 Since $$ \bar H^{i+1}\iso E^{\frac {i+3}2,\frac {i-1}2}, $$
(a)  has been proved.
\bigskip

For (b) we note that $l>1$ implies $g\ge r>r-l+1$ so we may use
[IV], Theorem 5.3.

The proof of (c) is similar to the proof of (a). We just mention
that, for all $i$ under consideration, we have the ``southwest''
isomorphisms
$$ \bar H^i=E^{i,0}\iso
E^{i-1,1}\iso\dots\iso E^{\frac {i-l}2+1,\frac {i+l}2-1},\ \bar
H^{i+1}=E^{i+1,0}\iso E^{i,1}\iso\dots\iso E^{\frac {i-l}2+2,\frac
{i+l}2-1}.
$$
\end{proof}

If $t$ is a negative integer, we obtain a similar result. We
indicate it without the proof which is completely analogous with
the proof of the previous theorem.
\begin{theorem}\label{FUNDAMENTALnegativ}Let $t<0$ be an integer.
Assume that $ 1\le r\le g$, and use the notation from above.
\begin{enumerate}
\item[\rm{(a)} ]Suppose that $t+l>0$ $($this implies $l>1$ $)$.
Then
\begin{enumerate}
\item[\rm{(i)}] $\bar{H}^{i}=0$ for $0\le i< \min(h,\max(2,t+l))$;
\item[\rm{(ii)}] if $2\le t+l<h$, then
$\bar{H}^{t+l}=H_\psi^{0,t+l-1}$; \item[\rm{(iii)}] if $i-t-l$ is
odd, and $t+l+1\le i < h-1$, then one has a natural exact sequence
$$
 0\to \bar{H}^i\to S_{\frac{i-t-l-1}2}(H^*)\otimes
S_{\frac{i+t+l-1}2} (C) \to
H_\psi^{\frac{i-t-l+1}2,\frac{i+t+l-1}2} \to \bar{H}^{i+1} \to 0.
$$
\end{enumerate}
\item[\rm{(b)}] Suppose that $t+l\le 0$. Then $\bar{H}^i= 0$ for
$i=0,\ldots,\min(2,h-1)$. For $i$ odd, $3\le i < \min(h-1,2r)$,
one has a natural exact sequence
$$ 0\to \bar{H}^i\to S_{\frac{i-1}2-t-l}(H^*)\otimes
S_{\frac{i-1}2} (C) \to H_\psi^{\frac{i+1}2,\frac{i-1}2} \to
\bar{H}^{i+1} \to 0. $$
\end{enumerate}
\end{theorem}

We supply Theorem \ref{FUNDAMENTAL} by some simple results
concerning the homology of $\mathcal{N}(t)$ at $N^h$.
\begin{proposition}\label{HEAD} As in Theorem {\rm
\ref{FUNDAMENTAL}} assume that $1\le r\le g$ and set $C=\Cok\psi$.
Let $\mu=\min (h,2r+1)$ and $t\ge \frac{\mu}2-1$.
\begin{enumerate}
\item[\rm{(a)}] There is an exact sequence
$$
0\to E^{\mu-1,1}\to \bar{H}^\mu \to H_\phi^{\mu,0};
$$
in particular, if $\mu<3$, then there is an exact sequence
$$
0\to \bar{H}^\mu \to H_\phi^{\mu,0}.
$$
\item[\rm{(b)}] For $\mu\ge 3$ odd, there is an exact sequence
$$
0\to E^{\mu-1,1}\to D_{t-\frac{\mu-1}2}(H) \otimes
S_{\frac{\mu-1}2} (C) \to H_\psi^{\frac{\mu+1}2,\frac{\mu-1}2},
$$
and \item[\rm{(c)}] for $\mu\ge 3$ even, there is an exact
sequence
$$ 0\to \bar{H}^{\mu-1} \to D_{t-\frac{\mu-2}2}(H) \otimes S_{\frac{\mu-2}2}
(C) \to H_\psi^{\frac{\mu}2,\frac{\mu-2}2} \to \bar{H}^\mu \to
H_\phi^{\mu,0} . $$
\end{enumerate}
\end{proposition}
\begin{remark} One may easily deduce similar sequences in case $t<\frac{\mu}2-1$.
\end{remark}
\begin{proof}(a) This follows immediately from (4) and the fact
that $H^{p,0}_\psi=H^{p,1}_\psi=0$ for all $p$.
\smallskip

(b) In order to cover the case $\mu=2r+1$, we  modify the first
diagram in the proof of Theorem \ref{FUNDAMENTAL}: the diagram
$$ \CD && 0 && 0 && 0 && 0 \\ && @VVV @VVV @VVV @VVV \\
@>>>\Im\partial_\psi^{\frac{\mu-1}2,\frac{\mu-3}2}
@>>>\Im\partial_\psi^{\frac{\mu+1}2,\frac{\mu-3}2}
@>>>\Im\partial_\psi^{\frac{\mu+3}2,\frac{\mu-3}2}
@>>>\Im\partial_\psi^{\frac{\mu+5}2,\frac{\mu-3}2}\\ && @VVV @VVV @VVV @VVV \\
@>>>\Im\partial_\psi^{\frac{\mu-1}2,\frac{\mu-3}2}
@>>>\Ker\partial_\psi^{\frac{\mu+1}2,\frac{\mu-1}2}
@>>>\Ker\partial_\psi^{\frac{\mu+3}2,\frac{\mu-1}2}
@>>>\Ker\partial_\psi^{\frac{\mu+5}2,\frac{\mu-1}2}\\ && @VVV @VVV
@VVV @VVV
\\ && 0@>>> H_\psi^{\frac{\mu+1}2,\frac{\mu-1}2}
@>>>H_\psi^{\frac{\mu+3}2,\frac{\mu-1}2}
 @>>> H_\psi^{\frac{\mu+5}2,\frac{\mu-1}2}\\
  && &&  @VVV @VVV @VVV
&& &&\\ && &&  0 && 0 && 0 &&
\endCD
$$
has exact columns. Then, as in the proof of Theorem
\ref{FUNDAMENTAL}, we obtain an exact sequence
\begin{equation*}
 0\to E^{\frac{\mu+1}2,\frac{\mu-1}2}\to D_{t-\frac{\mu-1}2}(H)\otimes
S_{\frac{\mu-1}2} (C) \to
\Ker(H_\psi^{\frac{\mu+1}2,\frac{\mu-1}2} \to
H_\psi^{\frac{\mu+3}2,\frac{\mu-1}2} ) \to
E^{\frac{\mu+3}2,\frac{\mu-1}2}.
\end{equation*}
 Since $
E^{\frac{\mu+1}2,\frac{\mu-1}2}\iso E^{\mu-1,1}$ in this case, we
are done. \smallskip

(c) Set $i=\mu-1$ and use the sequence $(*)$ in the proof of
Theorem \ref{FUNDAMENTAL} to get the exact sequence $$ 0\to
\bar{H}^{\mu-1} \to D_{t-\frac{\mu-2}2}(H)\otimes S_{\frac
{\mu-2}2}(C)\to H_\psi^{\frac{\mu}2,\frac{\mu}2}\to
E^{\frac{\mu+2}2,\frac{\mu-2}2}\to 0.$$  Since
$E^{\frac{\mu+2}2,\frac{\mu-2}2}\iso E^{\mu-1,1}, $ we can glue
this sequence and the sequence obtained under (a) to get the
result.
\end{proof}

\section{Length Formulas}

In this section {\it $R$ is a noetherian ring, and $M$ an
$R$-module which has a presentation
$$
\CD 0@>>>\mathcal{F}@>\chi>>\mathcal{G}@>>>M@>>>0
\endCD
$$
where $\mathcal{F}$, $\mathcal{G}$ are free modules of ranks $m$
and $n$.} Then in particular $r=n-m\ge 0$.
\smallskip

Let $\bar\lambda:M\to\mathcal{H}$ be an $R$-homomorphism into a
finite free $R$-module $\mathcal{H}$ of rank $l\le n$. By
$\lambda:\mathcal G\to\mathcal H$ we denote the corresponding
lifted map. In case $\grade I_\chi$ has the greatest possible
value $r+1$, Theorem 6.10 in [IV] provides a comparably
satisfactory description of the homology of the generalized Koszul
complex $\mathcal C_{\bar{\lambda}}(t)$
\begin{equation*}
0\to \bigwedge^{n} M \otimes S_{p}(\mathcal{H})^*
 \stackrel{\partial_{\bar{\lambda}}}{\to}\cdots\stackrel{\partial_{\bar{\lambda}}}{\to}
   \bigwedge^{t+l} M
\stackrel{\nu_{\bar{\lambda}}}{\to} \bigwedge^{t}
M\stackrel{\partial_{\bar{\lambda}}}{\to} \cdots
\stackrel{\partial_{\bar{\lambda}}}{\to} S_{t}(\mathcal{H}) \to 0,
\end{equation*}
depending on the grade of $I_\lambda$.

In the following we deal with another extremal case. We assume
that $\grade I_\chi=\dim R$. Then $R/I_\chi$ has finite length
$\ell(R/I_\chi)$. (Generally the length of an $R$-module $N$ is
denoted by $\ell(N)$).

First we dualize the sequence $\mathcal F\overset{\chi}\to\mathcal
G\overset{\lambda}\to \mathcal H$. Then we set $F=\mathcal F^*,\
G=\mathcal G^*,\ H=\mathcal H^*$ and $\psi=\chi^*,\
\phi=\lambda^*$, to attain the setup which we studied in the
previous sections.

We recall to some notation and facts from section 6 in [IV]. As
there we consider the upper part of the Koszul bicomplex $\mathcal
K_{.,.}(t)$ which we rewrite as

$$
\CD
 &&\vdots &&&&\vdots && \vdots &&\\
 &&@VVV  &&@VVV @VVV  \\
\cdots@>>> B_t^{0,-1}@>>>
\cdots @>d_\phi>>B_t^{t,-1} @>\pm\nu^{\phi}>>B_t^{t+1,-1}@>d_\phi>>\cdots\\
&&@VVV && @V\partial_\psi VV  @V\partial_\psi VV   \\
\cdots@>>>B_t^{0,0}@>>> \cdots@>>> B_t^{t,0} @>\pm\nu^{\phi}>>
B_t^{t+1,0}@>>>\cdots
\endCD
$$
where $$
 B^{0,0}_t=\begin{cases}
D_t(H)\otimes\bigwedge^{m}G\otimes S_0(F)^*&\quad\text{ if $0\le
t$},\\
S_0(H^*)\otimes\bigwedge^{t+l+m}G\otimes S_0(F)^*&\quad\text{ if
$-l\le t< 0$},\\
S_{-t-l}(H^*)\otimes\bigwedge^{m}G\otimes S_0(F)^*&\quad\text{ if
$t< -l$.}
\end{cases}$$
Set $M^p=\Coker (B_t^{p,-1}\overset{\partial_\psi}\to B_t^{p,0})$.
The canonical surjection $B_t^{p,0}\to M^p$ yields a complex
morphism
$$
\CD \cdots @>>> B_t^{-1,0}@>>> B_t^{0,0}@>>>\cdots @>d_\phi
>>B_t^{p,0}@>>>
B_t^{p+1,0}@>>>\cdots\\
 && @VVV @VVV &&
@VVV  @VVV\\
\cdots @>>>M^{-1}@>>>M^0 @>>>\cdots
@>\bar{d_\phi}>>M^p @>>> M^{p+1}@>>>\cdots\\
\endCD
$$ where the maps $\bar {d_\phi}$ are induced by $d_\phi$. The lower row is denoted by
$\mathcal M(t)$. By [IV], Proposition 6.8, there is an isomorphism
$\mathcal M(\rho-t)\to\mathcal C_{\bar{\lambda}}(t) $ whose
inverse composed with the induced complex map $\mathcal
M(\rho-t)\to \mathcal N(\rho-t)$ yields a complex map
$\mu:\mathcal C_{\bar\lambda}(t)\to \mathcal N(\rho-t)$.

\begin{theorem}\label{SubMaximalCaseHomology} Let $M$, $\bar\lambda$, and $\lambda$ be as above.
Set $\rho=r-l$, $g=\grade I_\chi,\ h=\grade I_\lambda$. Suppose
that $g=\dim R=r$. Equip $\mathcal C_{\bar{\lambda}}(t)$ with the
graduation induced by the complex morphism $\mu:\mathcal
C_{\bar{\lambda}}(t)\to\mathcal{N}(\rho-t)$. Then the homology
modules $\tilde H^{i}$ of $C^{^.}_{\bar{\lambda}}(t)$ have finite
length for $i\le \min(h-1,2r)$.
\smallskip

Set $C=\Coker\chi^*,\ S_0(C)=R/I_\chi$, and assume that $h>0$.

\begin{enumerate}
\item[\rm{(a)}] Let $l=1$. Then for all $t\in \mathbb{Z}$ and $i$
odd, $0<i<\min(h-1,2r)$,
$$
\ell(\tilde H^{i})-\ell(\tilde
H^{i+1})=\ell(S_{\frac{i-1}2}(C))-\ell(S_{\frac{i+1}2}(C)).
$$
\item[\rm{(b)}] Let $l>1$. We distinguish four cases.
\begin{enumerate}
\item[\rm{(i)}] For all $t\le \frac{\rho}2$ and $i$ odd,
$0<i<h-1$,
$$
 \ell(\tilde H^{i})-\ell(\tilde H^{i+1})=\begin{pmatrix}
r-t-\frac{i+1}2\\
l-1
\end{pmatrix}\ell(S_{\frac{i-1}2}(C))-\begin{pmatrix}
r-t-\frac{i+3}2\\
l-1
\end{pmatrix}\ell(S_{\frac{i+1}2}(C)).
$$
\item[\rm{(ii)}] Suppose that $\frac{\rho}2<t\le \rho$. If $i$ is
odd, $0<i<\min(h-1,2(\rho-t))$, then
$$
 \ell(\tilde H^{i})-\ell(\tilde H^{i+1})=\begin{pmatrix}
r-t-\frac{i+1}2\\
l-1
\end{pmatrix}\ell(S_{\frac{i-1}2}(C))-\begin{pmatrix}
r-t-\frac{i+3}2\\
l-1
\end{pmatrix}\ell(S_{\frac{i+1}2}(C)).
$$
If $i-l$ is even, $2(\rho-t)+l+2\le i<h-1$, then
\begin{align*}
 \ell(\tilde H^{i})-\ell(\tilde H^{i+1})  =\begin{pmatrix}
\frac{i+l}2-\rho+t-2\\
l-1
\end{pmatrix} & \ell(S_{\frac{i+l}2-1}(C))-\\&\begin{pmatrix}
\frac{i+l}2-\rho+t-1\\
l-1
\end{pmatrix}\ell(S_{\frac{i+l}2}(C)).
\end{align*}
If $2(\rho-t)+l+1< h$, then $\ell(\tilde
H^{2(\rho-t)+l+1})=\ell(S_{r-t}(C))$. Moreover
$$
\tilde H^{i}=
\begin{cases}
 S_{\rho-t}(C) &\quad\text{if}\ \ i=2(\rho-t)+1<h,\\
\quad 0 &\quad\text{if}\ \ 2(\rho-t)+2\le i<\min(h,2(\rho-t)+l+1).
\end{cases}
$$
\item[\rm{(iii)}] Suppose that $\rho<t<r$. If $i+r-t$ is odd,
$r-t+1\le i<h-1$, then
\begin{align*}
 \ell(\tilde H^{i})-\ell(\tilde H^{i+1})  =\begin{pmatrix}
\frac{i-r+t-3}2+l\\
l-1
\end{pmatrix} & \ell(S_{\frac{i+r-t-1}2}(C))-\\&\begin{pmatrix}
\frac{i-r+t-1}2+l\\
l-1
\end{pmatrix}\ell(S_{\frac{i+r-t+1}2}(C)).
\end{align*}
If $r-t< h$, then $\ell(\tilde H^{r-t})=\ell(S_{r-t}(C))$.
Moreover $\tilde H^{i}=0$ if $0\le i<\min(h,r-t)$.
\item[\rm{(iv)}] Suppose that $r\le t$ and $i$ odd, $0<i<h-1$.
Then
$$
 \ell(\tilde H^{i})-\ell(\tilde H^{i+1})=\begin{pmatrix}
t-\rho+\frac{i-3}2\\
l-1
\end{pmatrix}\ell(S_{\frac{i-1}2}(C))-\begin{pmatrix}
t-\rho+\frac{i-1}2\\
l-1
\end{pmatrix}\ell(S_{\frac{i+1}2}(C)).
$$
\end{enumerate}
\end{enumerate}
\end{theorem}
\begin{remark} Observe that in the above formulas we use the fact that
$h$ can reach its maximal value only if it is even. More
precisely, if $l\ge 2$ and $g=h=r$, then $l=2$ and $r$, $h$, $g$
are even (see Corollary 6.2 in [IV]). In this case $i$ odd,
$i<h-1$, means $i\le r-3$.

We further notice that for $h<\infty$ the formulas under (b) cover
the case in which $l=1$. We specified them for the readers
convenience since we shall apply the $l=1$ case in section 4.
\end{remark}

\begin{proof} If $r=0$, then $M=0$ since $\chi$ is injective.
So we may assume that $r\ge 1$.

The graduation induced by $\mathcal{N}(\rho-t)$ on $\mathcal
C_{\bar{\lambda}}(t)$ is completely determined by
$$C^0_{\bar{\lambda}}(t)=\begin{cases}
 \bigwedge^r M\otimes S_{\rho-t}(\mathcal H)^* &\quad \text{if
$t\le\rho$},\\ \bigwedge^t M \otimes S_0(\mathcal H)&\quad \text{if $\rho<t<r$},\\
 \bigwedge^r M \otimes S_{t-r}(\mathcal H) &\quad \text{if
$r\le t$}.\\ \end{cases}$$
  For $q>r$ the support of $\bigwedge ^qM$ is contained in the variety of $I_\chi$.
  Consequently $C^{i}_{\bar{\lambda}}(t)$ has finite length if $i<0$, which in turn implies that
$\tilde H^{i}$ has finite length. In particular, there remains
nothing to prove if $h=0$. Let $h>0$. Then $\rho\ge 0$ by
Proposition 5.1. in [IV].

By [IV], Proposition 6.9, $\mu:\mathcal
C_{\bar{\lambda}}(t)\to\mathcal{N}(\rho-t)$ induces the following
commutative diagram
$$
\begin{CD}
&& && 0 &&  0 &&  0\\
&& && @VVV  @VVV @VVV\\
C^{-1}_{\bar{\lambda}}@>\partial^{-1}_{\bar{\lambda}}>>C^{0}_{\bar{\lambda}}
@>\partial^{0}_{\bar{\lambda}}>>C^{1}_{\bar{\lambda}}@>>> C^{2}_{\bar{\lambda}}@>>> C^{3}_{\bar{\lambda}}\\
&& @V\mu_0 VV  @V\mu_{1} VV @V\mu_{2}VV @V\mu_{3}VV \\
0 @>>> N^0@>\bar d^0_\phi>>N^1 @>>> N^2@>>> N^3\\
&& @VVV  @VVV @VVV @VVV\\
0@>>> \Coker\mu_0@>\alpha>>\Coker\mu_1  @>>>  0@>>>  0\\
&& @VVV @VVV\\
&& 0 && 0
\end{CD}
$$
with exact columns.

For arbitrary $t$ and arbitrary $i$ the maps $\mu_i$ are
isomorphisms at all prime ideals which do not contain $I_\chi$.
Consequently $\Ker\mu_i$ and $\Coker\mu_i$ have finite length. In
particular $\Ker \mu_0$ equals the torsion submodule of
$C^{0}_{\bar{\lambda}}$ since $N^0$ is free. On the other hand,
$C^{1}_{\bar{\lambda}}$ is a torsion free module. So the torsion
submodule of $C^{0}_{\bar{\lambda}}$ is contained in $\Ker
\partial^{0}_{\bar{\lambda}}$. If we denote by $\bar \mu_0$ and $\bar
\partial^{0}_{\bar{\lambda}}$ the maps induced by $ \mu_0$ and $
\partial^{0}_{\bar{\lambda}}$ on $C^{0}_{\bar{\lambda}}/\Ker
\mu_0$, we get the following commutative diagram
$$
\begin{CD}
&& 0 && 0 &&  0 &&  0\\
&& @VVV @VVV  @VVV @VVV\\
&& C^{0}_{\bar{\lambda}}/\Ker \mu_0
@>\bar\partial^{0}_{\bar{\lambda}}>>C^{1}_{\bar{\lambda}}@>>> C^{2}_{\bar{\lambda}}@>>> C^{3}_{\bar{\lambda}}\\
&& @V\bar\mu_0 VV  @V\mu_{1} VV @V\mu_{2}VV @V\mu_{3}VV \\
0 @>>> N^0@>\bar d^0_\phi>>N^1 @>>> N^2@>>> N^3\\
&& @VVV  @VVV @VVV @VVV\\
0@>>> \Coker\mu_0@>\alpha>>\Coker\mu_1  @>>>  0@>>>  0\\
&& @VVV @VVV\\
&& 0 && 0
\end{CD}
$$
with exact columns.

Since $h>0$, $\bar d^0_\phi$ is injective. Then
$\bar\partial^{0}_{\bar{\lambda}}$ must be injective, which in
turn implies that $ \Ker \mu_0=\Ker \partial^{0}_{\bar{\lambda}}.
$ As $\tilde H^{0}$ is a factor of $\Ker
\partial^{0}_{\bar{\lambda}}$, we deduce that  $\tilde H^{0}$ has
finite length.

In case $\rho<t<r$, $\mu_0$ is injective by Proposition 6.9,(3) in
[IV]. So $\partial^{0}_{\bar{\lambda}}$ is injective and $\tilde
H^{0}=0$. In particular, we proved the only statement for $h=1$.

Assume that $h\ge 2$. If $h=\infty$, then $l=1$ (by Theorem 6.1 in
[IV]), and we refer to the next paragraph of the proof. So let
$\infty> h\ge 2$. Then, in particular, $r\ge 2$. The row
homologies at $N^0$ and $N^1$ vanish (see Theorems
\ref{FUNDAMENTAL} and \ref{FUNDAMENTALnegativ}). Of course,
$\tilde H^{1}=\Ker \alpha$ has finite length. We focus on the last
statements under (ii) and (iii) about $\tilde H^{1}$. Since
alternatively $t=\rho$, $t=r-1$, $\rho<t\le r-2$, another
application of Proposition 6.9,(3) in [IV] shows that
$\Coker\mu_1=0$ in all these cases. So
$$\tilde
H^{1}=\Coker\mu_0=
\begin{cases}
 S_{0}(C) &\quad\text{if}\ \ t=\rho,\\
 H^r(\mathcal C_{\psi}(1)) &\quad\text{if}\ \ t=r-1,\\
 0 &\quad\text{if}\ \ \rho<t\le r-2.
\end{cases}
$$
By Proposition 2.3 in [BV1] we have $\ell(H^r(\mathcal
C_{\psi}(1)))=\ell( S_{1}(C))$. With that we proved all claims for
$h=2$.

Now suppose that $h\ge 3$. If $l=1$, then the row homologies at
$N^0$, $N^1$ and $N^2$ vanish (see Theorems \ref{FUNDAMENTAL} and
\ref{FUNDAMENTALnegativ}), and we obtain
$$
\ell(\tilde H^{1})-\ell(\tilde H^{2})=\ell(\Ker
\alpha)-\ell(\Coker \alpha)=\ell(\Coker\mu_0)-\ell(\Coker\mu_1).
$$
But $\Coker\mu_0=S_{0}(C)$ and $\Coker\mu_1=H^r(\mathcal
C_{\psi}(1))$, so
$$
\ell(\tilde H^{1})-\ell(\tilde
H^{2})=\ell(S_{0}(C))-\ell(H^r(\mathcal C_{\psi}(1)).
$$
If $i$ is odd, $3\le i<h-1$, then we deduce directly from Theorem
\ref{FUNDAMENTAL},(a) and (c) and from Theorem
\ref{FUNDAMENTALnegativ},(b) that
$$
\ell(\tilde H^{i})-\ell(\tilde
H^{i+1})=\ell(S_{\frac{i-1}2}(C))-\ell(H^r(\mathcal
C_{\psi}(\frac{i+1}2)).
$$
Proposition 2.3 in [BV1] implies that $\ell(H^r(\mathcal
C_{\psi}(k)))=\ell( S_{k}(C))$ whenever $0\le k\le r$. It remains
to prove that $\tilde H^{h-1}$ has finite length if $h$ is even.
But Proposition \ref{HEAD},(c) provides an injection of $\tilde
H^{h-1}$ into a module of finite length. So we settled the case in
which $l=1$.

Let ($h\ge 3$ and) $l>1$. Only the case in which $\rho\le t<r$,
deserves special attention. The other cases are similar to the
case $l=1$. We computed already $\tilde H^0$ and $\tilde H^1$. If
$\rho\le t<r-2$ or $\rho=t=r-2$, then the row homologies at $N^0$,
$N^1$, $N^2$ vanish as $\Coker\mu_1$ does. So $\tilde
H^{2}=\Coker\mu_1=0$ in these cases. It remains to show that
$\ell(\tilde H^2)= \ell(S_2(C))$ if $\rho<t=r-2$. Of course, the
row homologies at $N^0$ and $N^1$ vanish. Since $\Coker\mu_1=0$,
we get that $\tilde H^{2}$ equals the row homology at $N^2$ which
is $H_\psi^{0,1}(-l+1)$ (see Theorem \ref{FUNDAMENTALnegativ}, (a)
(ii)). By Proposition 2.3 in [BV1] we obtain the desired length
equality. The remaining claims follow easily if one uses as
pattern the proof for the $l=1$ case.
\end{proof}

\smallskip
Let $\bar\lambda$ be as in the theorem, and denote by
$\tilde{\mathcal C_{\bar{\lambda}}}(t)$ the complex obtained from
$\mathcal C_{\bar{\lambda}}(t)$ by replacing
$C^i_{\bar{\lambda}}(t)$ with $0$ whenever $i<0$. As for $\mathcal
C_{\bar\lambda}(t)$ let $\tilde H^k$ denote the homology of
$\tilde{\mathcal C_{\bar{\lambda}}}(t)$ at $\tilde
{C^k_{\bar{\lambda}}}(t)$. Actually it differs from the
corresponding homology in $\mathcal C_{\bar\lambda}(t)$ only for
$k\le 0$.

\begin{corollary}\label{length formulas}
We adopt the assumptions and the notation of the first paragraph
in Theorem {\rm \ref{SubMaximalCaseHomology}}. As there we set
$C=\Coker\chi^*$ and $S_0(C)=R/I_\chi$.
\begin{enumerate}
\item[\rm{(a)}]If $h=\infty$, then
$$
\ell(S_0(C))=\ell(S_1(C))=\ldots=\ell(S_r(C)).
$$
\item[\rm{(b)}]If $h$ is odd and $t\le \frac{\rho}2$, then
$$
\sum_{k=0}^{h-1}(-1)^{k}\ell(\tilde H^k)=\begin{pmatrix}
r-t-\frac{h+1}2\\
l-1
\end{pmatrix}\ell(S_{\frac{h-1}2}(C)).
$$
\item[\rm{(c)}]If $h$ is even and $t\le \frac{\rho}2$, then
$$
\sum_{k=0}^{h-2}(-1)^{k}\ell(\tilde H^k)=\begin{pmatrix}
r-t-\frac{h}2\\
l-1
\end{pmatrix}\ell(S_{\frac{h-2}2}(C)).
$$
\end{enumerate}
In case $t> \frac{\rho}2$ one can easily deduce formulas similar
to {\rm (b)} and {\rm (c)}.\end{corollary}

\begin{proof}(a) If $h=\infty$, then $l$ must be $1$ (see Theorem 6.1 in [IV]).
We notice that this result is also an easy consequence of
Proposition 2.8 in [BV1].

(b) and (c) We may obviously suppose that $h>0$. So $\rho\ge 0$.
We have to prove that
$$
\ell(\tilde H^0(t))=\ell(D_{\rho-t}(H)\otimes S_0(C))
$$
if $t\le \frac{\rho}2$. From the proof of Theorem
\ref{SubMaximalCaseHomology} we deduce that
$$
\tilde H^0(t)\iso D_{\rho-t}(H)\otimes H^r(C^{^.}_{\psi}(0)),
$$
and from [BV1], Proposition 2.3 we draw that
$$
\ell(H^r(C^{^.}_{\psi}(0)))=\ell(S_0(C)).
$$
\end{proof}

\section{An application to quasi-homogeneous icis}

In this final section we extend the formulas of the previous
section in a very special case. Following the line of
argumentation in [BV1] we shall give a purely algebraic proof of
an old theorem due to Greuel (see [G], Proposition 2.5). Though
the length formulas and - at least implicitly - the proof of
Greuel's Theorem are contained in [HM], we present our
considerations as a byproduct of a more general approach.

We specialize to the case in which $R$ is a quasi-homogeneous
complete intersection with isolated singularity. More precisely,
we let $S=k[[X_1,\ldots,X_n]]$ where $k$ is a field of
characteristic zero, assign positive degrees $a_i$ to the
variables $X_i$, and set
$R=S/(p_1,\ldots,p_m)=k[[x_1,\ldots,x_n]]$ where the $p_i\in
(X_1,\ldots,X_n)^2$ form a regular sequence of homogeneous
polynomials of degrees $b_i$. By the Euler formula
$$
b_jp_j=\sum_{i=1}^n a_i\frac{\partial p_j}{\partial X_i}X_i.
$$
Since $\sum_{j=1}^mS\, p_j=\sum_{j=1}^mS\, (b_jp_j)$, the $b_jp_j$
may be viewed as defining elements for $R$. If we set
$p'_j=b_jp_j$ and $X'_i=a_iX_i$, we get
\begin{equation*}
p'_j=\sum_{i=1}^n\frac{\partial p_j}{\partial X_i}X'_i.\tag{$*$}
\end{equation*}
We suppose $m<n$ and $R_{\frak p}$ to be regular for all prime
ideals ${\frak p}$ different from the maximal ideal. As usual we
denote by $\Omega_{R/k}$ the module of K\"{a}hler-differentials of
$R$ over $k$. There is a presentation
$$
\CD 0@>>>\mathcal{F}@>\chi>>\mathcal{G}@>>>\Omega_{R/k}@>>>0
\endCD
$$
where $\mathcal{F}$, $\mathcal{G}$ are free $R$-modules of ranks
$m,\, n$ and $\grade I_\chi=r$. Moreover the Euler derivation
$\bar\lambda$ gives rise to an exact sequence
\begin{equation*}
\bigwedge^r\Omega_{R/k}\to\bigwedge^{r-1}\Omega_{R/k}\to\cdots\to\Omega_{R/k}\stackrel{\bar\lambda}{\to}
R\to k\to 0\tag{$**$}
\end{equation*}
which is in fact the non-negative grade part of $\mathcal
C_{\bar{\lambda}}$. Let $\lambda:\mathcal{G}\to R$ be the
corresponding lifted map. Set $\phi=\lambda^*$, $\psi=\chi^*$ as
above. As in the proof of Theorem \ref{SubMaximalCaseHomology} we
can complement ($**$) to an exact sequence
$$
0\to\tau(\bigwedge^r\Omega_{R/k})\to\bigwedge^r\Omega_{R/k}\to\bigwedge^{r-1}\Omega_{R/k}\to\cdots\to\Omega_{R/k}\to
R\to k\to 0
$$
where $\tau$ denotes the torsion submodule.

\begin{theorem}\label{Greuel isomorphism}Set $C=\Coker\psi$ and $S_0(C)=R/I_\chi$.
If $0\le i\le r-1$, then
$$
H^r(\mathcal C_{\psi}(i+1))\iso S_i(C).
$$
\end{theorem}

\begin{proof} We choose bases $g_1,\ldots,g_n$ for $G$ and $f_1,\ldots,f_m$ for
$F$ such that $\psi$ is represented by the matrix $(\frac{\partial
p_j}{\partial x_i})_{i,j}$, while $\phi$ is represented by
$(x'_1,\ldots,x'_n)$ (we denote by $\frac{\partial p_j}{\partial
x_i}$ the image in $R$ of $\frac{\partial p_j}{\partial X_i}$ and
by $x'_i$ the image of $X'_i$). Then we associate the bicomplex
$\mathcal{K}_{.,.}(0)$ with the sequence $R\overset\phi\to
G\overset\psi\to F$ (see section 1). The commutative diagram
$$
\begin{CD}
\cdots@>>>\bigwedge^{n-r} G
@>>>\cdots @>d_\phi>>\bigwedge^n G\iso R@>>>0\\
&& @V\pm\nu_0 VV && @V\pm\nu_r VV\\
0@>>>\bigwedge^0 G\iso R
@>>>\cdots @>d_\phi>>\bigwedge^r G@>>>\cdots\\
\end{CD}
$$
is the `middle' part of this bicomplex. First we prove the
following
\smallskip

\noindent{\bf Claim (1)}: {\it $\nu_r(g_1\wedge\ldots\wedge g_n)$
generates the homology in the second row at $\bigwedge^r G$.}
\smallskip

\noindent (For the original proof see the second part of the proof
of Theorem 3.1 in [HM]). We have
\begin{align*}
\nu_r(g_1\wedge\ldots\wedge g_n) & =g_1\wedge\ldots\wedge
g_n\leftharpoonup
\psi^*(f^*_1)\wedge\ldots\wedge\psi^*(f^*_m)\\
& =\sum_{\sigma}\epsilon(\sigma) \det_{1\le i,j\le m}
(\psi^*(f^*_j)(g_{\sigma(i)}))g_{\sigma(m+1)}\wedge\ldots\wedge
g_{\sigma(n)},
\end{align*}
where $\sigma$ runs through the set of permutations of $n$
elements which are increasing on the intervals $[1,m]$ and
$[m+1,n]$ (see [IV], section 2).

On the other hand there is a (non-canonical) complex isomorphism
$$
\begin{CD}
0@>>>\bigwedge^{0} G
@>>>\cdots @>d_\phi>>\bigwedge^n G@>>>0\\
&& @V\pm\Omega_0 VV && @V\pm\Omega_n VV\\
0@>>>\bigwedge^n \mathcal{G}
@>>>\cdots @>\partial_\lambda>>\bigwedge^0 \mathcal{G}@>>>0\\
\end{CD}
$$
induced by the isomorphism $\Omega_n:\bigwedge^n G\to R$,
$\Omega(g_1\wedge\ldots\wedge g_n)=1$. The lower row is the Koszul
complex associated with $\lambda$. If we denote by $H_i(R)$ the
row homology at $\bigwedge^i \mathcal{G}$, then
$H_m(R)\iso\bigwedge^m H_{1}(R)$ by a theorem of Tate and Assmus
(see Theorem 2.3.11 in [BH]). The relations ($*$) imply that
$H_{1}(R)$ is generated by the homology classes of the cycles
$\psi^*(f^*_j)$, $j=1,\ldots,m$ (see Chapter 2.3 in [BH]), so
$H_m(R)$ is generated by
$\psi^*(f^*_1)\wedge\ldots\wedge\psi^*(f^*_m)$. An easy
computation shows that
\begin{align*}
\Omega_r(\nu_r(g_1\wedge\ldots\wedge g_n))& =\pm\sum_{\sigma}
\det_{1\le i,j\le m}
(\psi^*(f^*_j)(g_{\sigma(i)}))g^*_{\sigma(1)}\wedge\ldots\wedge
g^*_{\sigma(m)}\\
&=\pm \psi^*(f^*_1)\wedge\ldots\wedge\psi^*(f^*_m),
\end{align*}
where $\sigma$ runs as above. Consequently
$\nu_r(g_1\wedge\ldots\wedge g_n)$ generates the homology at
$\bigwedge^r G$ as we claimed.

The complex map $\nu$ in the first diagram induces the commutative
diagram (D)
$$
\begin{CD}
&& 0 && 0 &&  0 &&\\ && @VVV @VVV  @VVV \\
0@>>>\bigwedge^r\Omega_{R/k}/\tau(\bigwedge^{r}\Omega_{R/k})
@>>>\bigwedge^{r-1}\Omega_{R/k}@>>>\bigwedge^{r-2}\Omega_{R/k}@>>>\cdots
\\ && @V\mu_0 VV
@V\mu_{1} VV @V\mu_2 VV\\ 0 @>>> N^0=R@>>>N^1@>\bar d_\phi>>N^2 @>>>\cdots\\
&& @VVV  @VVV @VVV\\ 0@>>>S_0(C) @>\alpha>>H^r(\mathcal
C_{\psi}(1)) @>>> 0\\ && @VVV @VVV\\ && 0 && 0
\end{CD}
$$
with exact columns where the $N^i$ are defined as in section 2. It
is the specialized version of the second diagram in the proof of
Theorem \ref{SubMaximalCaseHomology}. The homology of the lowest
row is denoted by $ h^{^.}$. The claim implies that
$\nu_r(\bigwedge^n G)\not\subset \Im d_\phi$. So
$\mu_r(R)\not\subset\Im{\bar d_\phi}$.

Let $r=1$. $R$ being a complete intersection, the homology of
$\mathcal D_\phi$ at $\bigwedge^r G$ is $k$ (for arbitrary $r\ge
1$). From the diagram
$$
\begin{CD}
&& 0 && 0\\
&& @VVV @VVV\\
0 @>>> N^0 @>>> N^1 @>>> 0\\
&& \parallel && @VVV @VVV\\
0 @>>> R @>>> C^{1,0} @>>> C^{2,0}\\
&& @VVV @VVV @VVV\\
&& 0 @>>> \Im\partial_\psi @>>> \Im\partial_\psi\\
&& && @VVV @VVV\\
&& && 0 && 0
\end{CD}
$$
with exact columns, we deduce that the homology of $\mathcal{N}$
at $N^1$ is also $k$, since the first non-trivial maps in the
second and in the third row are injective. So we obtain an exact
sequence
$$
\CD 0@>>> h^0@>>>k@>\beta>>k@>>> h^1@>>>0.
\endCD
$$
Since $\mu_r(R)\not\subset\Im \bar d_\phi$, $\beta$ must be an
isomorphism, and consequently $\alpha$ is an isomorphism.

Let $r=2$. First we show that $H^2(\mathcal C_{\psi}(1))\iso
S_0(C)$. The row homologies at $N^0$, $\Omega_{R/k}$, and $N^1$
vanish. Therefore $h^0=0$, and we get an exact sequence
$$
\CD 0@>>>h^1@>>>k@>\gamma>>H(N^2)
\endCD
$$
where $H(N^2)$ denotes the row homology at $N^2$. Because $\gamma$
is induced by $\mu_{2}$ and $\mu_r(R)\not\subset\Im \bar d_\phi$,
$\gamma$ must be injective, so $h^1=0$, and $\alpha$ is an
isomorphism. Next we show that $H^2(\mathcal C_{\psi}(2))\iso
S_{1}(C)$. Set $\mathfrak m=(x_1,\ldots,x_n)$. By Proposition 2.3
in [BV1] and the local duality theorem (see 3.5.8 in [BH]) we have
$$
H^2(C^{^.}_{\psi}(2))\iso\Ext^{1}(\bigwedge^{2}\Omega_{R/k},R)\iso(H^{1}_{\mathfrak
 m}(\bigwedge^{2}\Omega_{R/k}))^{\vee} \iso (S_0(C))^{\vee}\iso $$
$$
(H^2(C^{^.}_{\psi}(1)))^{\vee}\iso (H^{1}_{\mathfrak
 m}(\Omega_{R/k}))^{\vee}\iso\Ext^{1}(\Omega_{R/k},R)\iso
S_{1}(C).
$$

Finally let $r\ge 3$. Adopting the notation of section 2, we shall
prove
\smallskip

\noindent{\bf Claim (2)}: {\it If $1\le j\le \dfrac {r+1}2$, then
$E^{i,j}=0$ for $i=j,\ldots,r-j+1$, and $H^r(\mathcal
C_\psi(j))=S_{j-1}(C)$.}
\smallskip

\noindent We argue by induction on $j$. Let $j=1$. In the diagram
(D) the row homologies at $N^0$, $N^1$, $N^2$ vanish, so $\alpha$
must be an isomorphism which proves that $H^r(\mathcal
C_\psi(1)\cong S_0(C)$. Furthermore $\mathcal{N}$ has homology
only at $N^r$, namely $k$. Consider the commutative diagram
$$
\CD &&0&&0&&0&&0\\ && @VVV  @VVV  @VVV @VVV \\ \cdots@>>>
N^{r-2}@>>>N^{r-1} @>\bar d_{\phi}>>N^r @>>> 0\\
 && @VVV @VVV
@V\iota VV @VVV
\\ \cdots @>>> C^{r-2,0}@>>> C^{r-1,0}@>d_\phi>>C^{r,0}@>>> C^{r+1,0}@>>>\cdots
\\ &&
@VVV  @VVV  @VVV @VVV \\
\cdots@>>>\Im \partial_\psi^{r-2,0}@>>> \Im
\partial_\psi^{r-1,0}@>>> \Im
\partial_\psi^{r,0} @>>> \Im \partial_\psi^{r+1,0}@>>>\cdots\\
 && @VVV  @VVV
@VVV @VVV \\ &&0&&0&&0&&0\\
\endCD
$$
with exact columns. The second row is the Koszul complex
associated with $\phi$. It follows immediately that $E^{i,1}=0$
for $i=1,\ldots,r-2$. Furthermore we get an exact sequence
$$
\CD 0@>>>E^{r-1,1}@>>>k@>\beta>>k@>>>E^{r,1}@>>>0.
\endCD
$$
Since $\iota(N^r)=\nu_r(\bigwedge^n G)$ ($\mu_r$ is an
isomorphism), $\beta$ must be an isomorphism. We deduce
$E^{r-1,1}=E^{r,1}=0$.

Now let $1\le j+1\le\frac {r+1}2$. The commutative diagram
\scriptsize$$
\begin{CD}
&&0&&0&&0&&&&0\\ && @VVV  @VVV  @VVV &&@VVV \\
0@>>> \Im\partial_\psi^{j,j-1}@>>> \Im\partial_\psi^{j+1,j-1}@>>>
\Im\partial_\psi^{j+2,j-1}@>>>\ldots@>>>
\Im\partial_\psi^{r-j+1,j-1}&\ldots\\&& @VVV  @VVV  @VVV&& @VVV \\
0@>>>C^{j,j}@>>>C^{j+1,j}@>>>C^{j+2,j}@>>>\ldots@>>>C^{r-j+1,j}&\ldots\\
&& @VVV  @VVV  @VVV &&@VVV
\\0@>>>S_j(C)@>\gamma>>\Cok\partial_\psi^{j+1,j-1}@>\delta>>\Im\partial_\psi^{j+2,j}@>>>
\ldots@>>>\Im\partial_\psi^{r-j+1,j}&\ldots\\&& @VVV  @VVV @VVV&&
@VVV \\&&0&&0&&0&&&&0
\end{CD}
$$\normalsize
has exact columns and its second row is the Koszul complex of
$\phi$ tensored with $S_j(F)$, and therefore its homology vanishes
at least up to $C^{r-1,j}$. By our induction hypothesis the
homology of the first row vanishes up to
$\Im\partial_\psi^{r-j+1,j-1}$. In particular the row homologies
at $S_j(C)$ and $\Cok\partial_\psi^{j+1,j-1}$ vanish. Since $\Im
\partial_\psi^{j+1,j}\overset{\bar d_\phi}\to \Im
\partial_\psi^{j+2,j}$ is injective and is induced by $\delta$, we
obtain
\begin{equation*}
H^r(\mathcal C_\psi(j+1))\cong\Ker\delta=\Im\gamma\cong
S_j(C)\quad\text{and}\quad \Im\delta=\Im \bar d_\phi.
\end{equation*}
Together with the last equation we deduce $E^{i,j+1}=0$ for
$i=j+1,\ldots,r-j$. This proves our claim.

The claim implies that $H^r(\mathcal C_\psi(i+1))=S_{i}(C)$ for
$0\le i<\frac r2$. For $\frac{r}2\le i<r$, we use Proposition 2.3
in [BV1] and the local duality theorem as in the $r=2$ case to get
the desired isomorphism.
\end{proof}

\begin{corollary}\label{symetric length}With notation and the assumptions from above we obtain
$$\ell(S_0(C))=\ell(S_1(C))=\ldots=\ell(S_{r}(C)).$$
\end{corollary}

\begin{proof}Use the isomorphisms of Theorem \ref{Greuel isomorphism} and
Corollary 2.2 in [BV1].
\end{proof}

\begin{remark}In [BV1], section 3, the length formulas of Corollary \ref{symetric length} have been
 proved for $r$ odd. The reader may also find a proof in [HM], Proposition 4.9.
Our approach  follows the line of [BV1]. The isomorphisms of
Theorem \ref{Greuel isomorphism} were previously obtained only for
$0\le i\le r-2$, and consequently only the formula
$\ell(S_0(C))=\ldots=\ell(S_{r-1}(C)).$
\end{remark}

\begin{proposition}\label{length-sequence}
Set $M_\phi=\Coker\phi$, and let $\bar{\psi}:M_\phi\to F$ be the
map induced by $\psi$. Then for the homology of $\mathcal
C_{\bar{\psi}}(r)$, the following holds:
$$
H^{0}(\mathcal C_{\bar{\psi}}(r))=H^{1}(\mathcal
C_{\bar{\psi}}(r))=0,\quad  H^{r+1}(\mathcal C_{\bar{\psi}}(r))=
S_r(C).
$$
If $i+r$ is odd, $0\le i\le r-1$, then
$$ \ell({H}^i(\mathcal C_{\bar{\psi}}(r)))=\ell({H}^{i+1}(\mathcal C_{\bar{\psi}}(r))). $$
If $r$ is even, then $ H^{2}(\mathcal C_{\bar{\psi}}(r))=0$. (To
avoid misunderstandings: the graduation of $\mathcal
C_{\bar{\psi}}(r)$ is fixed in such a way that position $0$ is
held by $\bigwedge^n M_{\phi}$.)
\end{proposition}

\begin{proof}The bicomplex
 $$ \CD
 &&0&&0&&0&&\\
 && @VVV @VVV @VVV  \\
 @>>> N^{r-1}
 @>>> N^{r}@>>>k@>>>0\\
 && @VVV @VVV @VVV \\
 @>>>\bigwedge^{r-1}G\otimes S_{0}(F) @>>> \bigwedge^{r}G\otimes S_{0}(F)
 @>>>\bigwedge^{r}M_\phi\otimes S_{0}(F)@>>>0\\
  && @VVV @VVV @VVV \\
 @>>>\bigwedge^{r-2}G\otimes S_{1}(F) @>>> \bigwedge^{r-1}G\otimes S_{1}(F)
 @>>>\bigwedge^{r-1}M_\phi\otimes S_{1}(F)@>>>0\\
  && @VVV @VVV @VVV \\
 @>>>\bigwedge^{r-3}G\otimes S_{2}(F) @>>> \bigwedge^{r-2}G\otimes S_{2}(F)
 @>>>\bigwedge^{r-2}M_\phi\otimes S_{2}(F)@>>>0\\
 && @VVV @VVV @VVV \\
 &&\vdots && \vdots && \vdots
\endCD
$$
arises from $\tilde{\mathcal C}_{.,.}(r)$ by truncation at the
$(r+1)$th column and taking cokernels. Its last column is just
$\mathcal C_{\bar\psi}(r)$. (Observe that $\bigwedge^{n}M_\phi\iso
k$.) From Claim (1) in the proof of Theorem \ref{Greuel
isomorphism} we draw that $H^{0}(\mathcal C_{\bar{\psi}}(r))=0$.

To prove $H^1(\mathcal C_{\bar\psi}(r))=0$, let $Q$ be the total
ring of fractions of $R$. Consider the commutative diagram
$$
\begin{CD}
\bigwedge^{r}M_\phi@>\partial_{\bar\psi}>>\bigwedge^{r-1}M_\phi\otimes F\\
@VVV@VVV\\\bigwedge^{r}M_\phi\otimes Q@>\partial_{\bar\psi}\otimes
Q>>\bigwedge^{r-1}M_\phi\otimes F\otimes Q
\end{CD}
$$
with canonical vertical arrows. The kernel of the left vertical
arrow is the torsion of $\bigwedge^{r}M_\phi$. The right vertical
arrow and $\partial_{\bar\psi}\otimes Q$ are injective since
$\bigwedge^{r-1}M_\phi$ is torsionfree and $\mathcal
C_{\bar\psi}(r)\otimes Q$ is split exact. So
$\Ker(\bigwedge^{r}M_\phi\overset{\partial_{\bar\psi}}\to\bigwedge^{r-1}M_\phi)\otimes
F$ is the torsion submodule of $\bigwedge^{r}M_\phi$. Next we note
that this torsion equals $H^r(\mathcal D_\phi)$ where $\mathcal
D_\phi=\mathcal D_\phi(t)$ is the dual version of the Koszul
complex associated with $\phi$. Of course $H^r(\mathcal
D_\phi)=H_m(R)$. Furthermore $H_1(R)=k^m$. Since $R$ is a complete
intersection, the Koszul algebra $H_{.}(R)$ is isomorphic with the
exterior algebra of $H_1(R)$. Therefore $H_m(R)=k$.  Since
$\bigwedge^n M_\phi=k$ and $\nu_{\bar\psi}:\bigwedge^n M_\phi\to
\bigwedge^{r}M_\phi$ is injective, we get $H^1(\mathcal
C_{\bar\psi}(r))=0$.

The isomorphism ${H}^{r+1}(\mathcal C_{\bar{\psi}}(r))\iso S_{r}
(C)$ is obvious.

By the same method we used to prove Theorem \ref{FUNDAMENTAL}, we
obtain exact sequences
$$ 0\to {H}^i(\mathcal C_{\bar{\psi}}(r))\to S_{\frac{i+r-1}2} (C) \to
H^r(\mathcal C_\psi(\frac{i+r+1}2)) \to {H}^{i+1}(\mathcal
C_{\bar{\psi}}(r)) \to 0.
$$
for $i+r$ odd, $0\le i\le r-1$. Since
$$
S_{\frac{i+r-1}2} (C) \iso H^r(\mathcal C_\psi(\frac{i+r+1}2))
$$
by Theorem \ref{Greuel isomorphism}, the length formula follows.
Moreover $H^2(\mathcal C_{\bar{\psi}}(r))=0$ if $r$ is even.

\end{proof}

\begin{remark}\label{exact-sequence} We conjecture that actually
$$
H^{i}(\mathcal C_{\bar{\psi}}(r))=
\begin{cases}
 S_r(C) &\quad\text{if}\ \ i=r+1,\\
 0 &\quad\text{otherwise}.
\end{cases}
$$
\end{remark}

The conjecture is true if $m=\rank F=1$. {\it Proof}. In the case
under consideration we have an the exact sequence of complex
morphisms
$$ \mathcal C_{\bar\psi}(r)[1]\overset{\iota}\to \mathcal C_{\psi}(r)\overset{\pi}\to\mathcal
C_{\bar\psi}(r)\to 0$$ where $[\ ]$ means shift and
$$\pi_p=\begin{cases} \bigwedge^n (G\to M_\phi)\quad&\text{for $p=0$}\\
\bigwedge^{r-p+1}(G\to M_{\phi})\quad&\text{for
$p=1,\ldots,r+1$}.\end{cases}
$$ Obviously $\iota_p$ is injective for $p>0$. Since
$H^p(\mathcal C_\psi(r))=0$ for $p=0,\ldots,r-1$ and $H^0(\mathcal
C_{\bar\psi}(r))=H^1(\mathcal C_{\bar\psi}(r))=0$, we obtain
$H^p(\mathcal C_{\bar\psi}(r))=0$ for $p=0,\ldots,r$ as desired.

\end{document}